\documentclass[11pt]{amsart}
\usepackage{lineno,hyperref}
\usepackage{enumerate}
\usepackage{amsthm}
\usepackage{amsmath}
\usepackage{amssymb}
\newcommand{\N}{N}
\newcommand{\R}{\mathbb{R}}
\newcommand{\om}{\omega}
\newcommand{\RN}{\R\sp\N}
\newcommand{\nint}{\int_{\RN}}
\newcommand{\V}{V}
\newcommand{\F}{R}

\newcommand{\mm}{M}
\newcommand{\kk}{j}

\newcommand{\D}{D}
\newcommand{\Ce}{\sigma}
\newcommand{\co}{c_1}

\newcommand{\ra}{\to}
\newcommand{\var}{\varepsilon}
\renewcommand{\D}{\nabla}

\newtheorem{theorem}{Theorem}[]
\theoremstyle{definition}
\newtheorem{remark}{Remark}[]
\newtheorem{lemma}{Lemma}[]
\newtheorem{prop}{Proposition}[]

\begin{document}

\title[Standing-wave solutions to Klein-Gordon systems]{Energy estimated frequencies of standing-wave solutions to non-linear Klein-Gordon systems in higher dimension}

\thanks{%
This work has been supported by the \textsl{PDE Research Group} of School of Mathematical Sciences of the University of Nottingham Ningbo China and funded by the \textsl{FoSE New Researchers Grant}}
\author{Daniele Garrisi}
\email{daniele.garrisi@nottingham.edu.cn}
\address{Sir Peter Mansfield Building\\
School of Mathematical Sciences\\
University of Nottingham Ningbo China\\
199 Taikang East Road\\
315100, Ningbo, China}
\begin{abstract}
In this work a system of non-linear elliptic equations is considered,
where the non-linear term is the sum of a quadratic form and a Sobolev
sub-critical term. An extra assumption is introduced on the sub-critical term,
which is minimal among the ones which guarantee the existence of standing-waves
obtained by estimating frequencies of minimizing sequences with the energy
functional.
\end{abstract}
\maketitle
\section{Introduction}
\thispagestyle{empty}
\noindent
In this work, we show the existence of solutions to the
following system of \(m\) equations
\begin{linenomath}
\begin{equation}
\label{eq:elliptic-system}
-\Delta u_\kk + m_\kk\sp 2 u_\kk + 
\partial_\kk \F(u) = \om_\kk\sp 2 u_\kk,\quad 1\leq\kk\leq\mm.
\end{equation}
\end{linenomath}
We consider non-linearities which allow us to obtain solutions as minima
of the energy on the charge, as it has been done in \cite{BF09,BBBM10}. 
The energy functional is defined as
\begin{linenomath}
\begin{equation*}
E(\om,u) = \frac{1}{2} \sum_{\kk = 1}\sp\mm 
\nint\left(|\D u_\kk|\sp 2 + \om_\kk\sp 2 u_\kk\sp 2\right)dx 
+ \nint\V(u)dx,
\end{equation*}
\end{linenomath}
where
\begin{linenomath}
\begin{equation*}
\V(u) := \F(u) + \frac{1}{2}\sum_{\kk = 1}\sp\mm m_\kk \sp 2 u_\kk\sp 2.
\end{equation*}
\end{linenomath}
The constraints are defined as
\begin{linenomath}
\begin{equation*}
M_\Ce = \{(\om,u)\in \R\sp\mm\times H\sp 1(\R\sp\N,\R\sp\mm)\mid
C_\kk (\omega,u) = \Ce_\kk\}
\end{equation*}
\end{linenomath}
where
\begin{linenomath}
\begin{equation*}
C_\kk (\omega,u) = \om_\kk\nint |u_\kk|\sp 2dx.
\end{equation*}
\end{linenomath}
Existence is proved by showing that for some constraints \(\sigma\), Palais-Smale
sequences are bounded \(\R\sp\mm\times H^1 (\R\sp\N,\R\sp\mm)\), as in \cite{BF09}. Such boundedness follows from
\(\omega_\kk\leq m_{\kk}\), which holds if the energy of \((\omega,u)\)
is close enough to the infimum. In \cite{BF09}, this estimate is obtained
by selecting contraints \(M_\sigma\) such that in at least one point \((\omega,u)\) the hylenic charge, defined in \cite{BBBM10}
as 
\(\Lambda(\omega,u) := E(\omega,u)/C(\omega,u)\)
is strictly smaller than \(m\). A similar method has been applied in systems of equations, in \cite{Gar14}, where \(C\) is replaced by 
the sum of all the functions \(C_\kk\) for \(1\leq \kk\leq\mm\). In this work, we use a different strategy to dominate
\(\omega_\kk\). Using the assumption that \(V\geq 0\), which is also included in the references \cite{BF09,BBBM10} and \cite{Gar14}, 
one can show that 
\begin{linenomath}
\begin{equation}
\label{eq.domination}
\omega_\kk < \frac{2E(\omega,u)}{C_\kk(\omega,u)},\quad 1\leq\kk\leq\mm.
\end{equation}
\end{linenomath}
The following assumptions are required: 
\(\F\) is a real-valued continuously differentiable function on \(\mathbb{C}^m\) and
\(
\F(z) = \F(|z_1|,|z_2|,\dots,|z_m|)
\)
for every \(z\in\mathbb{C}^m\). We also assume that
\begin{linenomath}\begin{equation*}
\label{eq:A1}
\tag{$ A_1 $}
|\D\F(z)|\leq\co(|z|\sp{p - 1} + |z|\sp{q - 1}),\quad\F(0) = 0,\quad
2 < p\leq q < 2\sp *
\end{equation*}\end{linenomath}
\begin{linenomath}\begin{equation*}
\label{eq:A2}
\tag{$ A_2 $}
\V\geq 0.
\end{equation*}\end{linenomath}
When \(\mm\geq 2\), there exists \(u\) such that \(u_i > 0\) for every \(1\leq i\leq\mm\) and
\begin{linenomath}\begin{equation*}
\label{eq:A3}
\tag{$ A_3 $}
\V(u) - \frac{1}{8}\frac{\prod_{i = 1}^m m_i\sp 2 u_i\sp 2}%
{\sum_{|I|=m-1}\prod_I m_j\sp 2 u_j\sp 2} < 0.
\end{equation*}\end{linenomath}
When \(M=1\), we assume that there exists \(u > 0\) such that \(8\V(u) - m^2 u^2 < 0\).
To state our main result, we introduce the notation
\(
M_\sigma ^r := M_\sigma ^r\cap \big(\mathbb{R}^\mm _{++}\times H^1 _r (\mathbb{R}^N,\mathbb{R}^\mm)\big)
\),
where \(\mathbb{R}^\mm_{++} := \{u\in\mathbb{R}^\mm\mid u_i > 0\text{ for every } 1\leq i\leq M\}\).
\begin{theorem}
\label{thm:main}
Let $ \F $ fulfill assumptions (\ref{eq:A1}--\ref{eq:A3}). Then, there
exists $ L > 0 $ such that for every $ \Ce $ with $ |\sigma|\geq L $
\begin{enumerate}[(i).]
\item there holds \(2\inf_{M_\Ce^r}(E)\leq\min_{1\leq\kk\leq\mm} \{\Ce_\kk m_\kk \}\)
\item $ E $ satisfies the Palais-Smale condition at level $ c $ for every
\begin{linenomath}\begin{equation*}
c\in \Big[\inf_{M_\Ce^r}(E),\frac{1}{2}\min_{1\leq\kk\leq\mm}\{\Ce_\kk m_\kk \}\Big).
\end{equation*}\end{linenomath}
\end{enumerate}
\end{theorem}
This existence
result addresses the dimension \(N\geq 3\), as the proof relies on the
Radial Lemma, \cite{Str77,BL83a}.
\(A_1-A_3\) are tailored to the strategy used to prove
that if \((\omega_n,u_n)\) is a Palais-Smale sequence, then 
\(\liminf_{n\to\infty} \omega_n^\kk < m_\kk\), and are not meant
to be the sharpest assumptions for the existence of a solution to
the elliptic system \eqref{eq:elliptic-system}, especially 
for scalar equations, as pointed out in Remark~\ref{rk.1}. 
Critical points of \(E\) constrained to \(M_\sigma\) correspond to standing-wave
solutions to the system of non-linear Klein-Gordon equations
\begin{linenomath}
\begin{equation*}
\partial^2_t \phi(t,x) -\Delta \phi_\kk (t,x) + m_\kk\sp 2 \phi_\kk(t,x) + 
\partial_\kk \F(\phi(t,x)) = 0
\end{equation*}
\end{linenomath}
for \(1\leq\kk\leq\mm\). Systems of non-linear Klein-Gordon equations can use to describe models of relativistic
particle, where each particle interacts with the field and all the other particles. model interactions of each particle with
the field and the interaction with other particles, and their stability have been investigated in 
\cite{BBBM10}. When systems of two or more equations are considered, \(A_1-A_3\) are not able
to guarantee the stability of the ground-state. An example is provided in Remark~\ref{rk.2}, where
it is shown that minimizing sequences do not exihibit a concentration behaviour in the sense
of \cite{Lio84a}. For the stability of the ground-state to nonlinear
Klein-Gordon systems and Schr\"odinger systems, we refer to \cite{Gar12} and \cite{LNW16}.
\section{Properties of the functional $ E $}
In the following proposition we recall some well-known properties of the functional
$ E $.
\begin{prop}
\label{prop:properties}
Suppose that $ \F $ fullfils the assumption \eqref{eq:A1}. Then,
\begin{enumerate}[(i)]
\item $ E $ is continuously differentiable,
\item $ E $ is coercive,
\item for every bounded sequence $ (u_n) $ in \(H^1 _r\), up to extract a subsequence
\begin{linenomath}
\begin{equation*}
\nint (\nabla\F (u_{n}(x)) - \nabla\F(u_m(x)))\cdot(u_{n}(x) - u_m(x))dx\ra 0.
\end{equation*}
\end{linenomath}
\end{enumerate}
\end{prop}
\begin{proof}
(i). The continuity and differentiability of $ E $ follows from applications
and techniques of 
\cite[Theorem~2.2 and 2.6, p.\,16,~17]{AP93} for bounded domains. We 
do not include the proof as the reference applies
to unbounded domains and non-linearities satisfying \eqref{eq:A1}.

\noindent(ii).~Let \((\om,u)\) be a point of the constraint \(M_\Ce\)
and set \(E = E(\om,u)\). By \eqref{eq:A2}, we have
\begin{linenomath}
\begin{equation}
\label{eq:prop:properties-5}
E\geq\frac{1}{2}\sum_{j=1}^\mm \omega_j^2 \|u_j\|_{L^2}^2 = 
\frac{1}{2}\sum_{j=1}^\mm \omega_j\sigma_j,\quad
\|\D u\|_{L\sp 2} \sp 2\leq 2E.
\end{equation}
\end{linenomath}
Therefore,
\begin{linenomath}
\begin{equation}
\label{eq:prop:properties-1}
\om_\kk \leq\frac{2E}{\sigma_\kk}
\end{equation}
\end{linenomath}
for every \(1\leq\kk\leq \mm\). We set $ m := \min\{m_\kk\mid 1\leq\kk\leq\mm\} $. By \eqref{eq:A1} there exists 
$ \var > 0 $ such that
\begin{linenomath}
\begin{equation}
\label{eq:prop:properties-2}
\V(u)\geq \frac{m ^2 |u|^2}{4}
\end{equation}
\end{linenomath}
for every \(u\) such that \(|u|\leq\var\). We have
\begin{linenomath}\begin{equation*}
E\geq\int_{|u|\geq\var} \V(u)dx + \int_{|u| < \var} \V(u)dx.
\end{equation*}\end{linenomath}
From \eqref{eq:prop:properties-2}, it follows that
\begin{linenomath}
\begin{equation}
\label{eq:prop:properties-3}
\|u\|_{L^2 (|u| < \var)} ^2\leq \frac{4E}{m^2}.
\end{equation}
\end{linenomath}
On the other hand, by the Sobolev inequality, there exists \(k > 0\)
such that
\begin{linenomath}
\begin{equation}
\label{eq:prop:properties-4}
\begin{split}
\int_{|u|\geq\var} |u|\sp 2 dx &= \var\sp{2 - 2\sp *}\int_{|u|\geq\var}
\var\sp{2\sp * - 2} |u|\sp 2 dx\leq 
\var\sp{2 - 2\sp *} \int_{|u|\geq\var} |u|\sp{2\sp *}dx\\
&\leq k\sp{2\sp *} \var\sp{2 - 2\sp*}\|\D u\|_{L\sp 2} \sp{2\sp *}.
\end{split}
\end{equation}
\end{linenomath}
From \eqref{eq:prop:properties-3}, \eqref{eq:prop:properties-4}
and \eqref{eq:prop:properties-5} it follows that
\begin{linenomath}
\begin{equation*}
\|u\|\sp 2 _{L\sp 2}\leq 
\frac{4E}{m ^2} + 2 k\sp{2\sp *} \var\sp{2 - 2\sp*} E.
\end{equation*}
\end{linenomath}
Along with \eqref{eq:prop:properties-1}, we obtained that the
sub-levels of $ E $ are bounded subsets in \(\mathbb{R}^m\times H^1 _r (\mathbb{R}^n)\).
Therefore, $ E $ is coercive.

\noindent(iii). Using the Rellich-Kondrakov Theorem, from \cite[Radial Lemma]{Str77} it follows
that the inclusion \(H^1 _r (\mathbb{R}^N)\subset L^p _r(\mathbb{R}^N)\cap L^q _r(\mathbb{R}^N)\) is compact
for every \(2 < p < q < 2^*\). 
Therefore, up to extract a subsequence, we can suppose that there exists \(u\) in \(H^1 _r\)
such that \(u_n - u\) converges to zero in $ L\sp p (\mathbb{R}^N) \cap L\sp q (\mathbb{R}^N) $, where
\(p\) and \(q\) are as in \eqref{eq:A1}. By \eqref{eq:A1} and
the H\"older inequality, the following inequalities
\begin{linenomath}\begin{equation*}
\begin{split}
&\,\nint |(\nabla\F(u_n) - \nabla\F(u))(u_n - u)|dx
\leq\nint c_1(|u_n|\sp{p - 1} + |u|\sp{p - 1}) |u_n - u|dx
\\ 
+&\nint c_1 (|u_n|\sp{q - 1} + |u|\sp{q - 1}) |u_n - u|dx\leq
c_2 (\|u_n - u\|_{L\sp p} + \|u_n - u\|_{L\sp q})
\end{split}
\end{equation*}\end{linenomath}
hold for some $ c_2 = c_2(c_1,p,q,A) > 0 $, where
\begin{linenomath}\begin{equation*}
A = \sup_{n\geq 1}\max\big\{\|u_n\|^{p - 1}_{L^p},\|u_n\|^{q-1}_{L^q}\big\}.
\end{equation*}\end{linenomath}
Therefore, the rightend-side converges to zero.
\end{proof}
\section{The assumption \eqref{eq:A3}}
\label{sec.2}
In this section we explain the role of the assumption \eqref{eq:A3}
in obtaining estimates from above of \(\omega_\kk\) for every \(1\leq\kk\leq \mm\).
\begin{lemma}
\label{lem.minmax}
For every \(1\leq\kk\leq\mm\), let \(k_j\) be a positive real number and let
\(h_\kk (x) := \|x\|^2 - k_\kk x_\kk\) be the real-valued function defined on 
\(\mathbb{R}^\mm_{+}\). Then
the function 
\(
g(x) := \max_{1\leq \kk\leq\mm} h_\kk (x)
\)
achieves the infimum in the interior point of \(\mathbb{R}^\mm_{+}\) and 
\begin{linenomath}
\begin{equation}
\label{eq.minmax}
\min(g) = -\frac{\prod_{\kk=1}^\mm k_\kk^2}{4\sum_{|I|=\mm-1}\prod_{\kk\in I} k_\kk^2} = -\left(4\displaystyle\sum_{\kk=1}^\mm\frac{1}{k_\kk^2}\right)^{-1}.
\end{equation}
\end{linenomath}
\end{lemma}
\begin{proof}
Since \(h_\kk\) is a convex quadratic function, it is bounded below and coercive, \(g\) also is. Therefore,
the minimum is achieved at some point \(x^*\). Given \(x\) in \(\partial(\mathbb{R}^\mm_+)\),
there exists \(1\leq\kk\leq\mm\) such that \(x_\kk = 0\). Therefore, 
\(h_\kk(x)\geq 0\) implying \(\inf_{\partial(\mathbb{R}^\mm_+)} g\geq 0\).
To inspect the behaviour of \(g\) on the interior, we use the induction on \(\mm\). When \(\mm=1\) and \(x^*\) is an interior point,
we have \(0 = h_1'(x^*) = 2x^* - k_1x^*\), implying \(x^*=\frac{k_1}{2}\) and \(h_1(x^*)=-\frac{k_1^2}{4}\) proving
that the infimum is achieved in the interior of \(\mathbb{R}_+\). Suppose that \eqref{eq.minmax} holds for a given \(\mm\geq 1\). 
We study the case \(\mm + 1\). 
First of all, we show that there are \(1\leq a < b \leq\mm + 1\) such that
\(h_{a} (x^*) = h_b (x^*)\). On the contrary, there exists \(i_1\) such that
\(g(x^*) = h_{i_1} (x^*) > h_j (x^*)\) for every \(j\neq i_1\).
Then, there exists \(\delta > 0\) such that \(g(x) = h_{i_1}(x) > h_j(x)\) for every \(x\in B(x^*,\delta)\).
Since \(x^*\) is a minimum of \(g\), \(x^*\) is a local minimum of \(h_{i_1}\). Therefore, 
\begin{linenomath}\begin{equation*}
0 = \big(\nabla_x h_{i_1}(x^*)\big)_s = 2x_{i_1}^{*} - k_s\delta_{i_1 s}x_s,\quad 1\leq s\leq \mm + 1
\end{equation*}\end{linenomath}
implying \(x_j^* = 0\) for every \(j\neq i_1\) and \(x_{i_1}^* = \frac{k_{i_1}}{2}\). Thus \(h_{i_1}(x^*) = - \frac{k_{i_1}^2}{4} < 0 < h_j(x_*)\)
for every \(j\neq i_1\) contradicting the assumption that \(h_{i_1}(x^*) > h_j(x^*)\) for every \(j\neq i_1\). 
Let \(a < b\) be such that \(x_a^* k_a = x_b^* k_b\). 
For every \(x\) in 
\(
\Omega_{ab} := \{x\mid x_a k_a = x_b k_b\}
\)
there holds
\begin{linenomath}\begin{equation*}
\begin{split}
h_j (x) &= x_a^2  + x_b^2 + \sum_{s\neq a,b} x_s^2 - k_j x_j 
= \frac{k_a^2 x_a^2}{k_b^2} + x_a^2 + \sum_{s\neq a,b} x_s^2 - k_j x_j\\
&= \left(\frac{k_a^2}{k_b^2} + 1\right)x_a^2 + \sum_{s\neq a,b} x_s^2 - k_j x_j,\quad
\end{split}
\end{equation*}\end{linenomath}
Without loss of generality, we can suppose that \(a=1\) and \(b=2\). 
We define \(K\in\mathbb{R}^\mm_{++}\) and new variables \(X_\kk\) 
as
\begin{linenomath}\begin{equation*}
K_\kk = 
\begin{cases}
k_1\left(\frac{k_1^2}{k_2^2} + 1\right)^{-\frac{1}{2}} & \text{ if } \kk = 1\\
K_\kk = k_{\kk + 1} & \text{ if } 2\leq\kk\leq\mm
\end{cases}
\end{equation*}
\end{linenomath}
and
\begin{linenomath}
\begin{equation*} 
X_\kk = 
\begin{cases}
\left(\frac{k_1^2}{k_2^2} + 1\right)^{\frac{1}{2}}x_1  & \text{ if } \kk = 1\\
x_{\kk + 1} & \text{ if } 2\leq\kk\leq\mm.
\end{cases}
\end{equation*}\end{linenomath}
For every \(1\leq\kk\leq \mm\), we define \(H_\kk\) on \(\mathbb{R}^\mm_+\) as 
\(H_\kk (X) := \|X\|^2 - K_\kk X_\kk\).
Therefore,
\begin{linenomath}\begin{equation*}
\begin{split}
\min_{\mathbb{R}^{\mm+1}_+}\max_{1\leq\kk\leq\mm+1} h_\kk &= 
\min_{\mathbb{R}^{\mm+1}_{++}\cap{\Omega_{12}}}\max_{1\leq\kk\leq\mm+1} h_\kk\\
&=\min_{\mathbb{R}^{\mm}_+}\max_{1\leq\kk\leq\mm} H_\kk = -\left(4\sum_{\kk=1}^\mm\frac{1}{K_\kk^2}\right)^{-1}.
\end{split}
\end{equation*}\end{linenomath}
From
\begin{linenomath}\begin{equation*}
\begin{split}
\sum_{\kk=1}^\mm\frac{1}{K_\kk^2} &= \frac{1}{K_1^2} + \sum_{\kk = 2}^\mm\frac{1}{K_\kk^2} = 
\frac{1}{k_1^2}\left(\frac{k_1^2}{k_2^2} + 1\right) + \sum_{\kk = 2}^\mm\frac{1}{k_{\kk + 1}^2}\\
&= 
\frac{1}{k_1^2} + \frac{1}{k_2^2} + \sum_{\kk = 2}^\mm\frac{1}{k_{\kk + 1}^2} = 
\sum_{\kk = 1}^{\mm + 1}\frac{1}{k_\kk^2}
\end{split}
\end{equation*}\end{linenomath}
completing the induction step.
\end{proof}
\begin{prop}
Assumption \eqref{eq:A3} holds if and only if there exists \((\omega,u)\) in \(\mathbb{R}^{\mm}_{++}\times\mathbb{R}^{\mm}_{++}\) 
such that
\begin{equation}
\label{eq.A3-equivalent}
\tag{\(A_3'\)}
2\V(u) + \sum_{\kk = 1} ^{\mm} \om_\kk\sp 2 u_\kk\sp 2 - 
m_i \om_i u_i \sp 2 < 0,\text{ for all } 1\leq i\leq\mm.
\end{equation}
\end{prop}
\begin{proof}
Suppose that \eqref{eq.A3-equivalent} holds. 
We set 
\begin{linenomath}\begin{equation*}
f_\kk (\omega,u) := \sum_{\kk = 1} ^{\mm} \om_\kk\sp 2 u_\kk\sp 2 - 
m_i \om_i u_i \sp 2
\end{equation*}\end{linenomath}
and
\(g(\omega,u) := \max\{f_1 (\omega,u),f_2(\omega,u),\dots,f_\mm(\omega,u)\}\).
The condition \eqref{eq.A3-equivalent} can be restated as
\begin{equation}
\label{eq.A3-equivalent.2}
\inf_u\Big(2\V(u) + \inf_{\omega\in\mathbb{R}^\mm_+} g(\omega,u)\Big) < 0.
\end{equation}
Given \(u\in\mathbb{R}^\mm_{+}\), we consider the define constants \(k_\kk = m_\kk u_\kk\) and variable changes
\(x_\kk = \omega_\kk u_\kk\) for \(1\leq\kk\leq\mm\).
Therefore, \(f_\kk(\omega,u) = \|x\|^2 - k_\kk x_\kk =: h_\kk(x)\) for every \(1\leq\kk\leq\mm\). From Lemma~\ref{lem.minmax},
\begin{linenomath}\begin{equation*}
\inf_{\omega\in\mathbb{R}^\mm_+} g(\omega,u) = -\left(4\sum_{j=1}^\mm\frac{1}{k_j^2}\right)^{-1} = 
-\frac{1}{4}\frac{\prod_{\kk=1}^m m_\kk^2 u_\kk^2}{\sum_{|I|=m-1}\prod_{j\in I} m_j^2 u_j^2}.
\end{equation*}\end{linenomath}
Therefore, \eqref{eq.A3-equivalent.2} is equivalent to claim the existence of \(u\in\mathbb{R}^\mm_{++}\) such that
\eqref{eq:A3} holds.
\end{proof}
\section{Existence of stationary points on $ M_\sigma $}
\label{sec.3}
Hereafter, we assume that $ \Ce_\kk > 0 $ for every $ 1\leq\kk\leq\mm $.
\begin{prop}
\label{prop:palais-smale}
Let $ (\om_n,u_n)_{n = 1} \sp\infty\subset M_\sigma ^r $ be a Palais-Smale 
sequence such that $ \lim_{n\to\infty}\om_n\sp\kk < m_\kk $ for every 
\(1\leq\kk\leq\mm\).
Then $ (u_n)_{n = 1} \sp\infty $ has a converging subsequence.
\end{prop}
\begin{proof}
By (ii) of Proposition~\ref{prop:properties}, $ (u_n)_{n = 1} \sp\infty $
is bounded. Let \(L\) be such that \(\|u_n\|^2_{H^1}\leq L\) for every \(n\geq 1\).
Thus, we can suppose that $ u_n\rightharpoonup u $ in 
$ H\sp 1 _r (\RN) $.
Since $ (\om_n,u_n) $ is a Palais-Smale sequence, 
there are sequences \((\lambda_n \sp\kk)_{n = 1} \sp\infty\subset\R\)
and \((\eta_n,\xi_n)_{n = 1}\sp\infty\subset H^1_r \times\R\sp\mm\).
such that
\begin{linenomath}\begin{equation}
\label{eq:prop:palais-smale-1}
E'(\omega_n,u_n) = \sum_{\kk = 1}\sp\mm \lambda_n \sp\kk
C_\kk'(\omega_n,u_n) + (\eta_n,\xi_n),\quad (\eta_n,\xi_n)\ra 0.
\end{equation}\end{linenomath}
We apply both functionals in 
\eqref{eq:prop:palais-smale-1} to vectors in $ \R\sp\mm\times\{0\} $ and obtain
\(
\omega_n\sp\kk\|u_n\sp\kk\|_{L^2} ^2 = 
\lambda_n \sp\kk\|u_n\sp\kk\|_{L^2} ^2 + \eta_n \sp\kk
\),
whence \(\lambda_n \sp\kk  = \omega_n\sp\kk - \eta_n \sp\kk\omega_n \sp\kk/\Ce_\kk\).
Now we apply them to $ (0,\phi)\in \{0\}\times H\sp 1 _r $ and obtain
\begin{linenomath}\begin{equation*}
\begin{split}
\sum_{\kk=1}^\mm (\nabla u_n^j,\nabla\phi_j)_{L^2} &+ \sum_{\kk=1}^\mm
\nint \D\V(u_n)\phi dx + 
\sum_{\kk = 1}\sp\mm (\omega_n \sp\kk)^2 (u_n \sp\kk ,\phi_\kk)_{L^2}\\
&- 2\sum_{\kk = 1}\sp\mm\lambda_n \sp\kk  \omega_n \sp\kk
(u_n \sp\kk ,\phi_\kk)_{L^2} = \sum_{\kk=1}^\mm(\xi_n^\kk,\phi_j)_{L^2}.
\end{split}
\end{equation*}\end{linenomath}.
Substituting \(\lambda_n\sp\kk\) in the equality, we obtain
\begin{equation*}
\begin{split}
\nint\nabla u_n ^j\nabla\phi_j dx &+ \nint\partial_j\V(u_n)\phi_j dx - 
\sum_{\kk=1}^\mm(\om_n \sp\kk)\sp 2 (u_n \sp\kk,\phi_\kk)_{L^2} \\
&= \sum_{\kk=1}^\mm \left[(\xi_n\sp\kk,\phi_j)_{L^2}
- \frac{2\eta_n \sp\kk %
(\om_n \sp\kk)\sp 2}%
{\Ce_\kk} (u_n \sp\kk ,\phi_\kk)_{L^2}\right]
\end{split}
\end{equation*}
for every \(1\leq\kk\leq\mm\).
Because $ \omega_n $ converges to $ \omega $, from
\eqref{eq:prop:palais-smale-1} and the equality above it follows that
\begin{equation*}
\nint\nabla u_n ^j\nabla\phi_j dx + \nint\partial_j \V(u_n)\phi_j dx - 
\om_\kk\sp 2 (u_n \sp\kk,\phi_\kk)_{L^2} = \sum_{\kk = 1}^\mm (g_n^j,\phi_j)_{L^2}
\end{equation*}
where 
\begin{equation}
\label{eq:prop:palais-smale-4}
g_n^j = \xi_n^\kk - \frac{2\eta_n \sp\kk %
(\om_n \sp\kk)\sp 2}%
{\Ce_\kk} u_n \sp\kk = o(1)\text{ in } H^1 _r.
\end{equation}
We apply the equality above for \(n,m\in\mathbb{N}\) and 
\(\phi_j = u_n^\kk - u_m^\kk\), take the difference, and then the
sum over \(1\leq\kk\leq\mm\). Therefore,
\begin{equation*}
\begin{split}
&\|\D u_n^\kk - \D u_m^\kk\|_{L^2} ^2 + 
\sum_{\kk = 1}\sp\mm (m_\kk ^2 - \om_\kk ^2)
\|u_n\sp\kk - u_m\sp\kk\|_{L^2} ^2 \\
=& \sum_{\kk = 1}^m (g_n^\kk - g_m^\kk,u_n^\kk - u_m^\kk)_{L^2}
- \nint(\nabla\F(u_n) - \nabla\F(u_m))\cdot (u_n - u_m)dx.
\end{split}
\end{equation*}
Up to extract a subsequence, from (iii) of Proposition~\ref{prop:properties}
and from the assumption \(\omega_\kk < m_\kk\), the lefthand-side is a Cauchy
sequence. Therefore, a subsequence of \((u_n)_{n=1}\sp\infty\) converges in 
\(H^1 _r\).
\end{proof}
\begin{prop}
\label{prop:hylomorphy}
Suppose condition \eqref{eq:A3} holds. Then for every \(\alpha > 0\), there
exists \(\sigma\) such that \(\sigma_i \geq\alpha\) for every \(1\leq i\leq\mm\)
and $ 2\inf_{M_\Ce^r} E < \Ce_i m_i  $ for every $ 1\leq i\leq\mm $.
\end{prop}
\begin{proof}
Given $ (\om,u)\in\R\sp\mm_{++}\times\R\sp\mm_{++} $, we define
the \(H^1\) functions \(v_r\sp\kk(x) := u_\kk\chi_{\overline{B}(0,r)} + (1 + r - |x|)\chi_{\overline{C}(r,r+1)}\),
where \(C(r,r + 1) := B(0,r + 1)\cap\overline{B}^c (0,r)\) for every \(r > 0\), and \(\chi_A\) is the characteristic function of
a subset \(A\subset\mathbb{R}^N\).
By standard calculations, we have \(\D v_r\sp\kk = -\frac{u_\kk x}{|x|}\chi_{\overline{C}(r,r+1)}\)
and
\begin{linenomath}
\begin{equation*}
\|v_r\sp\kk\|_{L\sp 2}^2 = \mu(B_1) r\sp\N u_\kk\sp 2 + O(r\sp{\N - 1}),\quad
\nint \F(v_r)dx = \mu(B_1) r\sp\N \F(u) + O(r\sp{\N - 1})
\end{equation*}
\end{linenomath}
\begin{linenomath}
\begin{equation*}
\|\D u_r\sp\kk\|_{L\sp 2}^2 = O(r\sp{\N - 1}),
\end{equation*}
\end{linenomath}
where $ B_1 $ is the unit ball of $ \RN $ and \(\mu\) is the Lebesgue measure. Therefore,
\begin{linenomath}\begin{equation*}
\frac{2E(\om,v_r)}{C_i (\om,v_r)} = 
\frac{%
\sum_{\kk = 1}\sp\mm (m_\kk\sp 2 + \om_\kk\sp 2)
u_\kk\sp 2 + 2\F(u)}{\om_i u_i\sp 2} + o(1).
\end{equation*}\end{linenomath}
Taking the limit as $ r\ra +\infty $, the rightend-side becomes
\begin{equation}
\label{eq:prop:hylomorphy-1}
\frac{1}{\om_i u_i\sp 2}\cdot\left(%
2\F(u) + \sum_{\kk = 1} \sp\mm (m_\kk^2 + \om_\kk\sp 2) u_\kk\sp 2\right).
\end{equation}
From \eqref{eq.A3-equivalent} and the definition of \(V\), \eqref{eq:prop:hylomorphy-1} is smaller than $ m_i $.
Thus,
\begin{linenomath}\begin{equation*}
\frac{2E(\om,v_r)}{C_i (\om,v_r)} < m_i,\quad 1\leq i\leq\mm,
\end{equation*}\end{linenomath}
which completes the proof of the proposition.
\end{proof}
\begin{proof}[Proof of Theorem~\ref{thm:main}]
From Proposition~\ref{prop:hylomorphy}, there exists 
\(\sigma\) such that \(2\inf_{M_\Ce^r} E < \Ce_\kk m_\kk\)
for every \(1\leq\kk\leq\mm\).
Let \((\omega_n,u_n)_{n=1}^\infty\subseteq M_\sigma ^r\) be a minimizing sequence.
By the Ekeland Principle, \cite[Theorem~5.1,~p.~48]{Str90} there exists a Palais-Smale subsequence for $ E $ on $ M_\Ce^r $.
By (ii) of Proposition~\ref{prop:properties}, such sequence is bounded.
Then, we can suppose that $ \om_n\ra\om $. By inequality 
\eqref{eq:prop:properties-1}, we have
\begin{linenomath}\begin{equation*}
\om_n \sp\kk\leq\frac{2E(\omega_n,u_n)}{\Ce_\kk}.
\end{equation*}\end{linenomath}
Taking the limit as $ n\ra\infty $, we obtain $ \om_\kk < m_\kk $ for every $ 1\leq\kk\leq\mm $. 
By Proposition \ref{prop:palais-smale}, we can suppose that the sequence converges strongly to some \((\omega,u)\).
From (i) of Proposition~\ref{prop:properties},
\begin{linenomath}
\begin{equation*}
o(1) + \inf_{M_\sigma^r} (E) = E(\omega_n,u_n) = o(1) + E(\omega,u)
\end{equation*}
\end{linenomath}
implying \(E(\omega,u) = \inf_{M_\sigma^r}(E)\). By the Symmetric Criticality Principle, \cite[\S 0]{Pal79},
\((\omega,u)\) is a critical point of \(E\) constrained to \(M_\sigma\), thus a weak solution to the
elliptic system \eqref{eq:elliptic-system}.
\end{proof}
\begin{remark}
\label{rk.1}
Condition \eqref{eq:A3} is not a necessary for the existence of solutions to
the elliptic system \eqref{eq:elliptic-system}. In \cite{BBBM10} it has been proved that 
the existence of \(u > 0\) such that \(\F(u) < 0\) is sufficient, while the assumption \eqref{eq:A3}
reads \(\F(u) + \frac{3m^2 u^2}{8} < 0\).
\end{remark}
\begin{remark}
\label{rk.2}
Minima on \(M_\sigma\) may not exist. When \(\mm = 2\), and
\(N=3\), one
can just consider \(m_1 = m_2 = 1\) and define
\begin{linenomath}\begin{equation*}
\V(z_1,z_2) = \frac{|z|^2}{2} + \F(|z_1|) + \F(|z_2|),\quad \F(u) := \frac{u^2}{4(u + 1)^2} - \frac{u^2}{2}.
\end{equation*}\end{linenomath}
\(\V\) is smooth and satisfies \eqref{eq:A2} as \(V\geq 0\). \(\F\) satisfies \eqref{eq:A1} with \(p = 3\) and any \(q\) real number bigger than \(3\).
\eqref{eq:A3} holds as there exists \((u_1,u_2)\in\mathbb{R}^2 _{++}\), with \(u_1 = u_2\), such that
\begin{linenomath}\begin{equation*}
\frac{u_1^2}{4(u_1 + 1)^2} + \frac{u_2^2}{4(u_2 + 1)^2} - \frac{u_1^2 u_2^2}{8(u_1^2 + u_2^2)} < 0.
\end{equation*}\end{linenomath}
Minima of \(E\) over \(M_{(\Ce_1,\Ce_2)}\) exist. In fact, since
there are not coupling terms, one can deal with it using \cite[Theorem~2.8]{BBBM10}
for scalar equations. In fact \(W(z) := R(|z|) + \frac{|z|^2}{2}\) satisfies all the assumptions \cite[\(H_0-H_2\),~page~3]{BBBM10} and \cite[\(H_3\),~page~10]{BBBM10}: 
calculations show that \(W(0) = W'(0) = 0\) and \(W''(0) = 1\) implying
\(H_0\). Since 
\(
\inf\frac{W(s)}{2s^2} = 0,
\)
\(H_1\) holds, and \(W\geq 0\) implies \(H_2\). Finally, \(R\) is two-times
continuously differentiable and \(R''\) satifies the combined power-type
estimate in \(H_3\) with \(p = q = 3\). If \((u_1,u_2)\) is a minimum, \((u_1,u_2 (\cdot + y_n))\) is a minimizing sequence. One can
choose \((y_n)\) diverging, so the sequence does not concentrate. 
\end{remark}
\def\cprime{$'$} \def\cprime{$'$} \def\cprime{$'$} \def\cprime{$'$}
  \def\cprime{$'$} \def\cprime{$'$} \def\cprime{$'$} \def\cprime{$'$}
  \def\cprime{$'$} \def\polhk#1{\setbox0=\hbox{#1}{\ooalign{\hidewidth
  \lower1.5ex\hbox{`}\hidewidth\crcr\unhbox0}}} \def\cprime{$'$}
  \def\cprime{$'$} \def\cprime{$'$}
\providecommand{\bysame}{\leavevmode\hbox to3em{\hrulefill}\thinspace}
\providecommand{\MR}{\relax\ifhmode\unskip\space\fi MR }
\providecommand{\MRhref}[2]{%
  \href{http://www.ams.org/mathscinet-getitem?mr=#1}{#2}
}
\providecommand{\href}[2]{#2}

\thispagestyle{empty}
\end{document}